\documentclass[10pt,reqno]{amsart}

\usepackage{hyperref}
\usepackage{amssymb}
\usepackage{amsmath}
\usepackage{amsfonts}
\usepackage{graphicx,cite}
\usepackage{latexsym}
\usepackage{enumerate}
\usepackage{mathrsfs}
\usepackage{mathtools, mathdots}
\usepackage{booktabs}
\usepackage{xcolor}
\usepackage{caption}
\DeclarePairedDelimiter{\floor}{\lfloor}{\rfloor}
\usepackage{enumitem}
\setlist[itemize]{leftmargin=*}
\setlist[enumerate]{leftmargin=*}

\theoremstyle{plain}

\newtheorem{theorem}[equation]{Theorem}

\theoremstyle{definition}

\newtheorem{definition}{Definition} 
\newtheorem{example}[equation]{Example}

\newcommand{\0}{\textcolor{gray}{0}}

\newcommand{\C}{\mathbb{C}}
\newcommand{\F}{\mathcal{F}}
\newcommand{\Z}{\mathbb{Z}}
\newcommand{\M}{\mathbf{M}}

\newcommand{\K}{\mathcal{K}}
\newcommand{\X}{\mathcal{X}}

\newcommand{\T}{\mathsf{T}}
\newcommand{\Irr}{\operatorname{Irr}}

\newcommand{\Tau}{\mathscr{T}}
\newcommand{\inner}[1]{\left< #1 \right>}
\newcommand{\norm}[1]{\|#1\|}
\newcommand{\sgn}{\operatorname{sgn}}
\newcommand{\diag}{\operatorname{diag}}

\def\A{{\mathcal{A}}}
\renewcommand{\vec}[1]{{\boldsymbol{\bf #1}}}

\title[Supercharacters and the DFT, DCT, and DST]{Supercharacters and the discrete Fourier, cosine, and sine transforms}
\author{Stephan Ramon Garcia}
\address{\scriptsize Department of Mathematics, Pomona College, 610 N. College Ave., Claremont, CA 91711} 
\email{stephan.garcia@pomona.edu}
\urladdr{\url{http://pages.pomona.edu/~sg064747}}
\thanks{Partially supported by a David L. Hirsch III and Susan H. Hirsch Research Initiation Grant.
First author partially supported by National Science Foundation Grant DMS-1265973.}

\author{Samuel Yih}

\begin{document}

\let\oldenumerate=\enumerate
\def\enumerate{
\oldenumerate
\setlength{\itemsep}{5pt}
}

\let\olditemize=\itemize
\def\itemize{
\olditemize
\setlength{\itemsep}{5pt}
}


\maketitle
\begin{abstract}  
Using supercharacter theory, we identify the matrices that are diagonalized by the discrete cosine and discrete sine transforms, respectively. Our method affords a combinatorial interpretation for the matrix entries.
\end{abstract}

\section{Introduction}

The theory of supercharacters was introduced by P.~Diaconis and I.M.~Isaacs in 2008 \cite{SCT},
generalizing earlier seminal work of C.~Andr\'e \cite{An95, An01, An02}.  
The original aim of supercharacter theory was to provide new tools for handling the character theory of intractable groups,
such as the unipotent matrix groups $U_n(q)$.  Since then, supercharacters have appeared in the study of
combinatorial Hopf algebras \cite{Aguiar}, 
Schur rings \cite{Hendrickson,HuJo08}
and their combinatorial properties \cite{DiTh09, ThVe09, Th10},
and exponential sums from number theory \cite{RSS, GHM, SESUP}.

Supercharacter techniques permit us to identify the algebra of matrices that 
are diagonalized by the discrete Fourier transform (DFT) and discrete cosine transform (DCT),
respectively.  A natural modification handles the discrete sine transform (DST).
Although the matrices that are diagonalized by the DCT or DST have been studied previously
\cite{DPDCT,GMDST,SCP,ADCT}, we further this discussion in several ways. 

For the DCT, we produce a novel combinatorial description of the matrix entries
and obtain a basis for the algebra that has a simple combinatorial interpretation.
In addition to recapturing results presented from \cite{DPDCT}, we are
also able to treat the case in which the underlying cyclic group has odd order.

A similar approach for the DST runs into complications, but we can still 
characterize the diagonalized matrices by considering the ``orthocomplement'' of the DCT
supercharacter theory.
In special cases, the
diagonalized matrices are $\Tau$-class matrices \cite{SCP}, which first arose in the spectral theory of Toeplitz matrices
and have since garnered significant interest because of their computational advantages \cite{CHM,STBP, FDSM, FTTE}. 

For cyclic groups of even order, we recover results on \cite{SCP}. 
However, our approach also works if the underlying cyclic group has odd order.  This
is not as well studied as the even order case. In addition, we produce a second natural basis 
equipped with a novel combinatorial interpretation for the matrix entries. 

For all of our results, we provide explicit formulas for the matrix entries of the most general matrix
diagonalized by the DCT or DST, respectively.  

We hope that it will interest the supercharacter community to see that
its techniques are relevant to the study of matrix transforms that are traditionally the province of engineers,
computer scientists, and applied mathematicians.
Consequently, this paper contains a significant amount of exposition since we mean to bridge
a gap between communities that do not often interact. We thank the anonymous referee for 
suggesting several crucial improvements to our exposition.

\section{Preliminaries}

The main ingredients in this work are the theory of supercharacters
and the discrete Fourier transform (DFT), along with its offspring (the DCT and DST).  
In this section, we briefly survey some relevant definitions and ideas.

\subsection{Supercharacters}

The theory of supercharacters, which extends the classical character theory of finite groups,
was developed axiomatically by Diaconis--Isaacs \cite{SCT}, building upon earlier important work of
Andr\'e \cite{An95, An01, An02}.   It has since become an industry in and of itself.  We make no attempt to
conduct a proper survey of the literature on this topic.

\begin{definition}[Diaconis--Isaacs \cite{SCT}]
Let $G$ be a finite group, let $\X$ be a partition of the set $\Irr G$ of irreducible characters of $G$, 
and let $\K$ be a partition of $G$.  We call the ordered pair $(\X, \K)$ a \emph{supercharacter theory} if 
\begin{enumerate}
\item[(i)] $\K$ contains $\{0\}$, where $0$ denotes the identity element of $G$,
\item[(ii)] $|\X| = |\K|$,
\item[(iii)] For each $X \in \X$, the function  
$\sigma_X = \sum_{\chi \in X} \chi(0)\chi$ is constant on each $K \in \K$.
\end{enumerate}
The functions $\sigma_X$ are  \emph{supercharacters} and the elements $K$ of $\K$
 are  \emph{superclasses}.
\end{definition}

While introduced primarily to study the representation theory of non-abelian groups
whose classical character theory is largely intractable, recent work 
has revealed that it is profitable to apply supercharacter theory to the most elementary groups imaginable:
finite abelian groups \cite{RSS, GHM, SESUP, GNGP, GNSG, SAHS, HendricksonThesis, HendricksonNew}.

We outline the approach developed in \cite{SESUP}.  Although it is the ``one-dimensional'' case 
that interests us here, there is no harm in discussing things in more general terms.
Let $\zeta = \exp(-2\pi i /n)$, which is a primitive $n$th root of unity.
Classical character theory tells us that the set of irreducible characters of $G = (\Z/n\Z)^d$ is
\begin{equation*}
\Irr G = \{ \psi_{\vec{x}} : \vec{x} \in G \},
\end{equation*}
in which
\begin{equation*}
\psi_{\vec{x}}(\vec{y}) = \zeta^{\vec{x} \cdot \vec{y}}.
\end{equation*}
Here we write $$\vec{x} \cdot \vec{y} := \sum_{i=1}^d x_i y_i,$$ 
in which $\vec{x} = (x_1,x_2,\ldots,x_d)$ and $\vec{y}=(y_1,y_2,\ldots,y_d)$
are typical elements of $G$.  Since $\vec{x} \cdot \vec{y}$ is computed modulo $n$ it causes no 
ambiguity in the expression that defines $\psi_{\vec{x}}$. 
We henceforth identify the character $\psi_{\vec{x}}$ with $\vec{x}$.
Although this identification is not canonical (it depends upon the choice of $\zeta$),
this potential ambiguity disappears when we construct certain supercharacter theories on $G$.

Let $\Gamma$ be a subgroup of $GL_d(\Z/n\Z)$ that is closed under
the matrix transpose operation.  If $d = 1$, then $\Gamma$ can be any subgroup of the unit group $(\Z/n\Z)^{\times}$.
The action of $\Gamma$ partitions $G$ into $\Gamma$-orbits; we collect these orbits in the set
\begin{equation*}
\K = \{K_1,K_2,\ldots,K_N\}.
\end{equation*}
For $i = 1,2,\ldots,N$, we define
\begin{equation*}
\sigma_i := \sum_{\vec{x} \in K_i} \psi_{\vec{x}}.
\end{equation*}
The hypothesis that $\Gamma$ is closed under the transpose operation ensures that 
$\sigma_i$ is constant on each $K_i$ \cite[p.~154]{SESUP} (this condition is automatically
satisfied if $d = 1$).
For $i=1,2,\ldots,N$,
let $X_i = \{ \psi_{\vec{x}} : \vec{x} \in K_i\}$.  Then
\begin{equation*}
\X = \{X_1,X_2,\ldots,X_N\}
\end{equation*}
is a partition of $\Irr G$ 
and the pair $(\X,\K)$  is a supercharacter theory on $G$.

As an abuse of notation, we identify both the supercharacter and superclass partitions as $\{X_1,X_2,\ldots,X_N\}$
(such an identification is not always possible with general supercharacter theories). 
Since the value of each supercharacter $\sigma_i$ is constant on each superclass $X_j$,
we denote this common value by $\sigma_i(X_j)$.

Maintaining the preceding notation and conventions,
the following theorem links
supercharacter theory on certain abelian groups
and combinatorial-flavored matrix theory \cite[Thm.~2]{SESUP}.

\begin{theorem}[Brumbaugh, et.~al., \cite{SESUP}]\label{Theorem:Big}
For each fixed $z$ in $X_k$, let $c_{i,j,k}$ denote the number of solutions $(x_i,y_j) \in X_i \times X_j$ to $x+y = z$;
this is independent of the representative $z$ in $X_k$ which is chosen.
\begin{enumerate}
\item For $1\leq i,j,k,\ell \leq N$, we have
$$\displaystyle \sigma_i(X_{\ell}) \sigma_j(X_{\ell}) = \sum_{k=1}^N c_{i,j,k} \sigma_k(X_{\ell}).$$

\item The matrix 
\begin{equation}\label{eq:U}
U = \frac{1}{\sqrt{n^d}} \left[  \frac{   \sigma_i(X_j) \sqrt{  |X_j| }}{ \sqrt{|X_i|}} \right]_{i,j=1}^N
\end{equation}
is unitary ($U^* = U^{-1}$) and $U^4 = I$.

\item The matrices $T_1,T_2,\ldots,T_N$, whose entries are given by
\begin{equation}\label{eq:T}
[T_i]_{j,k} = \frac{ c_{i,j,k} \sqrt{ |X_k| } }{ \sqrt{ |X_j|} },
\end{equation}
each satisfy $T_i U = U D_i$, in which
$$D_i = \operatorname{diag}\big(\sigma_i(X_1), \sigma_i(X_2),\ldots, \sigma_i(X_N) \big).$$

\item Each $T_i$ is normal ($T_i^*T_i = T_i T_i^*$) and
the set $\{T_1,T_2,\ldots,T_N\}$ forms a basis for the algebra $\mathcal{A}$ of all $N \times N$ matrices 
$T$ such that $U^*TU$ is diagonal.
\end{enumerate}
\end{theorem}

The quantities $c_{ijk}$ are combinatorial in nature and are nonnegative integers that
relate the values of the supercharacters to each other.
Of greater interest to us is the unitary matrix $U$ defined in \eqref{eq:U}.  It is a normalized
``supercharacter table'' of sorts.  As in classical character theory, a suitable normalization of the
rows and columns of a character table yields a unitary matrix.  This suggests that
$U$ encodes an interesting ``transform'' of some type.
Theorem \ref{Theorem:Big} describes, in a combinatorial manner, 
the algebra of matrices that are diagonalized by $U$.  

This is the motivation for our work:
we can select $G$
and $\Gamma$ appropriately so that $U$ is either the discrete Fourier or discrete cosine
transform matrix.  Consequently, we can describe the algebra of matrices that 
are diagonalized by these transforms.
The discrete sine transform can be obtained as a sort of
``complement'' to the supercharacter theory corresponding to the DCT.  To our knowledge,
such complementary supercharacter theories have not yet been explored in the literature.

\subsection{The discrete Fourier transform}\label{Section:DFT}

It is hallmark of an important theory that even the simplest applications should
be of wide interest.  This occurs with the theory of supercharacters, for its most immediate byproduct is
the discrete Fourier transform (DFT), a staple in engineering and discrete mathematics.

A few words about the discrete Fourier transform are in order.
As before, let $G = \Z/n\Z$ and $\zeta= \exp(-2\pi i /n)$.
Let $L^2(G)$ denote the complex Hilbert space of all functions $f:G\to\C$, endowed with the inner product
\begin{equation*}
\inner{f,g} = \sum_{j=0}^{n-1} f(j) \overline{g(j)}.
\end{equation*}

The space $L^2(G)$ hosts two familiar orthonormal bases.
First of all, there is the \emph{standard basis} $\{\delta_0,\delta_1,\ldots,\delta_{n-1}\}$, which consists of the functions
\begin{equation*}
    \delta_j(k) = 
    \begin{cases}
        1 & \text{if $j=k$},\\
        0 & \text{if $j \neq k$}.
    \end{cases}
\end{equation*}
We work here modulo $n$, which explains our preference for the indices $0,1,\ldots,n-1$.
A second orthonormal basis of $L^2(G)$ is furnished by 
the \emph{exponential basis} $\{\epsilon_0,\epsilon_1,\ldots,\epsilon_{n-1}\}$, in which
\begin{equation*}
    \epsilon_j(\xi) = \frac{e^{2\pi i j \xi/n} }{ \sqrt{n}}.
\end{equation*}
The \emph{discrete Fourier transform} of $f \in L^2(G)$ is the function
$\widehat{f} \in L^2(G)$ defined by
\begin{equation*}
\widehat{f}(\xi) = \frac{1}{\sqrt{n}} \sum_{j=0}^{n-1} f(j) e^{-2\pi i j \xi / n} =  \inner{f, \epsilon_{\xi}}.
\end{equation*}
The choice of normalization varies from field to field.  We have selected the constant $1/\sqrt{n}$
so that the map $f \mapsto \widehat{f}$ is a unitary operator from $L^2(G)$ to itself.  
Indeed, the unitarity of the DFT follows from the fact that
\begin{equation*}
\widehat{\epsilon_j} = \delta_j, \qquad j=0,1,\ldots,n-1.
\end{equation*}
That is, the DFT is norm-preserving since it sends one orthonormal basis to another.
The matrix representation of the DFT with respect to the standard basis is
\begin{equation}\label{eq:FourierMatrix}
	F_n= 
    \frac{1}{\sqrt{n}} \begin{bmatrix} 
    1 & 1 & 1 & \cdots & 1 \\
    1 & \zeta & \zeta^2 & \cdots & \zeta^{n-1} \\
    1 & \zeta^2 & \zeta^4 & \cdots & \zeta^{2(n-1)} \\ 
    \vdots & \vdots & \vdots & \ddots & \vdots \\
    1 & \zeta^{n-1} & \zeta^{2(n-1)} & \cdots & \zeta^{(n-1)^2} \end{bmatrix}.
\end{equation}
This is the \emph{DFT matrix} of order $n$ (also called the \emph{Fourier matrix} of order $n$).

If we regard elements of $L^2(G)$ as column vectors, with respect to the standard basis,
then a short exercise with finite geometric series reveals that $\widehat{\delta_j} = \overline{\epsilon_j}$.
A little more work confirms that $F_n^2 = -I$ and hence $F_n^4 = I$.
Thus, the eigenvalues of $F_n$ are among $1,-1,i,-i$; the exact multiplicities
can be deduced from the evaluation of the quadratic Gauss sum, which is the trace
of $\sqrt{n} F_n$ \cite{Berndt}.

There are many compelling reasons why the discrete Fourier transform arises
in both pure and applied mathematics.  It would take us too far afield to go into details,
so we content ourselves with mentioning that the DFT arises in signal processing, number theory 
(e.g., arithmetic functions),
data compression, partial differential equations, and numerical analysis (e.g., fast integer multiplication).
A particularly fast implementation of the DFT, the \emph{fast Fourier transform} (FFT), was named
one of the Top 10 algorithms of the 20th century \cite{IEEE}.  
Although often credited to Cooley--Tukey (1965) \cite{Cooley}, 
the FFT was originally discovered by Gauss in 1805 \cite{Heideman}.
A valuable reference for all things Fourier-related is \cite{Kammler}.
The recent text of Stein and Shakarchi \cite{Stein} is a new classic
on the subject of Fourier analysis and it highly recommended for its
friendly and understandable approach.

How does the DFT relate to supercharacter theory?  Consider the following example,
which was first worked out in \cite{SESUP}.

\begin{example}[Discrete Fourier transform]
Let $G = \Z/n\Z$ and let $\Gamma = \{1\}$, the trivial subgroup of $(\Z/n\Z)^{\times}$, act upon $G$ by multiplication.
Then the $\Gamma$-orbits in $G$ are singletons:  $X_j = \{j\}$ for $j=0,1,2,\ldots,n-1$.
The corresponding supercharacters are classical exponential characters:
\begin{equation*}
\sigma_j(k) = \sum_{x \in X_j} \zeta^{x k} = \zeta^{j k}
\end{equation*}
and hence the unitary matrix $U$
from \eqref{eq:U} is the DFT matrix.  That is, $$U  = F_n.$$
Theorem \ref{Theorem:Big} permits us to identify the matrices
that are diagonalized by $U$.
With a little work, one can show that the matrices \eqref{eq:T} are
		\begin{equation*}
			[T_i]_{j,k} = 
			\begin{cases}
			0 & \text{if $k - j \neq i$},\\
			1 & \text{if $k - j = i$},
			\end{cases}
		\end{equation*}
		and they satisfy $T_i U= U D_i$, in which
		$$D_i = \operatorname{diag}(1, \zeta^{i}, \zeta^{2i},\ldots, \zeta^{(n-1)i}).$$ 
		The algebra $\mathcal{A}$ generated by the $T_i$ is the algebra of all
		$N \times N$ \emph{circulant matrices}
		\begin{equation*}
			\begin{bmatrix}
				c_0     & c_{N-1} & \cdots  & c_{2} & c_{1}  \\
				c_{1} & c_0    & c_{N-1} &         & c_{2}  \\
				\vdots  & c_{1}& c_0    & \ddots  & \vdots   \\
				c_{N-2}  &        & \ddots & \ddots  & c_{N-1}   \\
				c_{N-1}  & c_{N-2} & \cdots  & c_{1} & c_0 \\
			\end{bmatrix} .
		\end{equation*}
		More information about circulant matrices and their properties can be found in
		\cite[Sect.~12.5]{GarciaHorn}.
\end{example}

The preceding example shows that the discrete Fourier transform
arises as the simplest possible application of supercharacter theory.  If the action of
the trivial group $\{1\}$ on $\Z/n\Z$ already produces items of great interest,
it should be fruitful to consider actions of slightly-less trivial groups as well.
This motivates our exploration of the discrete cosine transform.

\section{Discrete cosine transform}\label{Section:DCT}

As before, we fix a positive integer $n$ and let $G = \Z/n\Z$.  Let
\begin{equation*}
    L^2_+(G) = \{ f \in L^2(G) : f(x) = f(-x)\quad \forall x\in G\}
\end{equation*}
and
\begin{equation*}
    L^2_-(G) = \{ f \in L^2(G) : f(x) = -f(-x)\quad \forall x\in G\}
\end{equation*}
denote the subspaces of even and odd functions in $L^2(G)$, respectively.  
Observe that $L^2_+(G)$ is invariant under the DFT, since, if $f$ is even, 
\begin{equation*}
    \widehat{f}(\xi) = \langle f, \epsilon_{\xi} \rangle = \frac{1}{\sqrt{n}} \sum_{j=0}^{n-1} f(j) e^{-2 \pi ij \xi/n} = \frac{1}{\sqrt{n}} \sum_{k=0}^{n-1} f(-k) e^{2 \pi ik \xi/n} = \langle f, \epsilon_{-\xi} \rangle = \widehat{f}(-\xi), 
\end{equation*}
and hence $\widehat{f}$ is even as well.  Since $L^2(G)$ is finite dimensional and the DFT is unitary,
it follows that $L^2_-(G) = L^2_+(G)^{\perp}$ is invariant under the DFT.  Consequently,
we have the orthogonal decomposition
\begin{equation*}
    L^2(G) = L^2_+(G) \oplus L^2_-(G),
\end{equation*}
in which both subspaces on the right-hand side are DFT-invariant.
The \emph{discrete cosine transform} (DCT) is the restriction of the DFT
to $L^2_+(G)$.  
Being the restrictions of a unitary operator (on a finite-dimensional Hilbert space)
to an invariant subspace, the DCT is a unitary operator on $L^2_+(G)$.
In a similar manner, 
the \emph{discrete sine transform} (DST) is the restriction of the DFT to $L^2_-(G)$.
It too is a unitary operator.

The DCT is a workhorse in engineering and software applications.
The MP3 file format, which contains compressed audio data, and the JPEG file format,
which contains compressed image data, make use of the DCT \cite{Gonzalez}.
These ``lossy'' file formats do not perfectly replicate the original source; that is,
some information is lost.  However, by judiciously eliminating high-frequency
components in the signal, one is able to produce sounds or images that
are, to human senses, virtually indistinguishable from the source.  Moreover,
this can be done in such a way that the final file size is much smaller than the original.

\begin{figure}
\centering
\includegraphics[width=\textwidth]{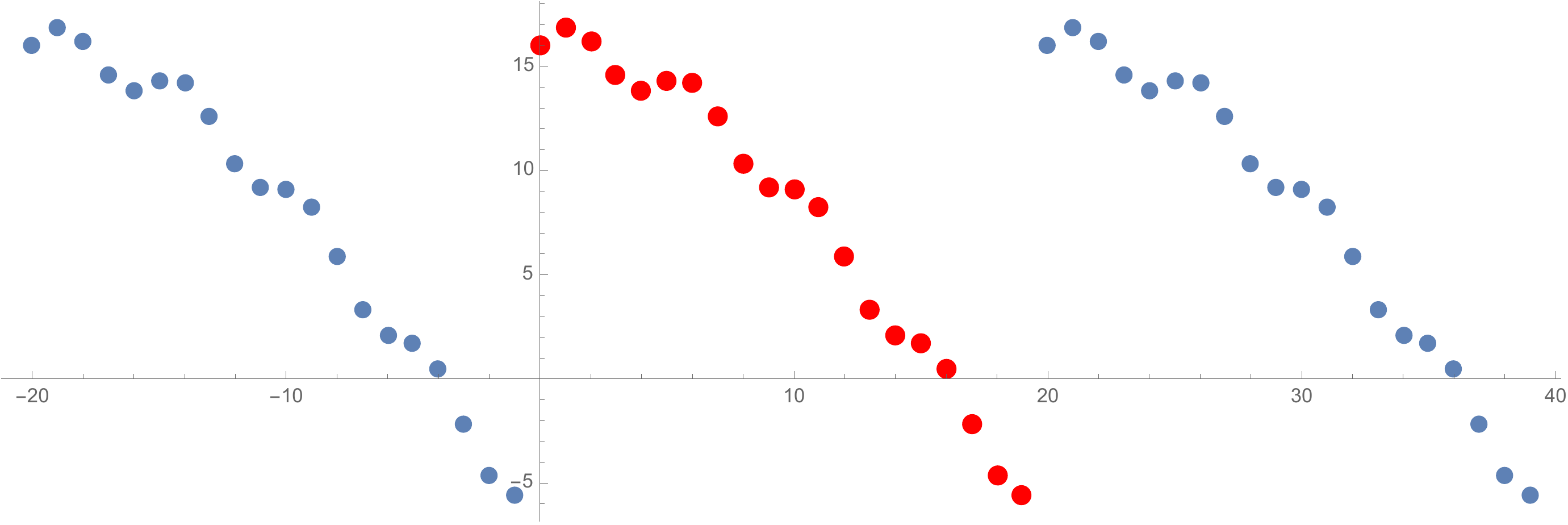}\\[15pt]
\includegraphics[width=\textwidth]{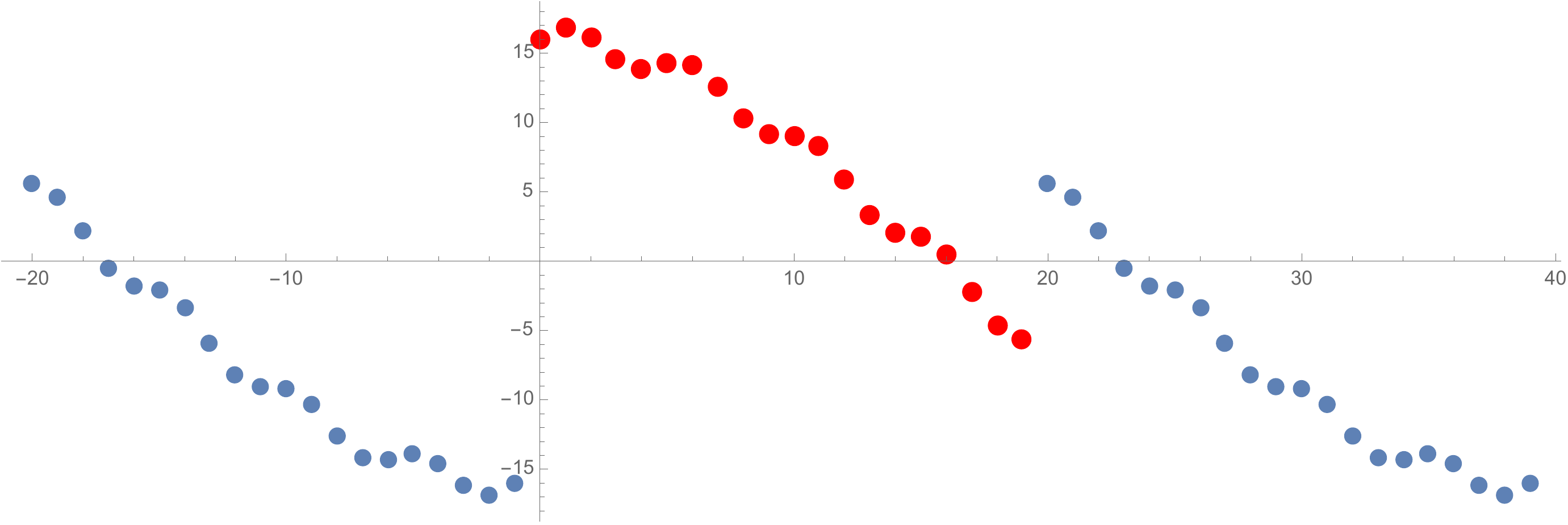}\\[15pt]
\includegraphics[width=\textwidth]{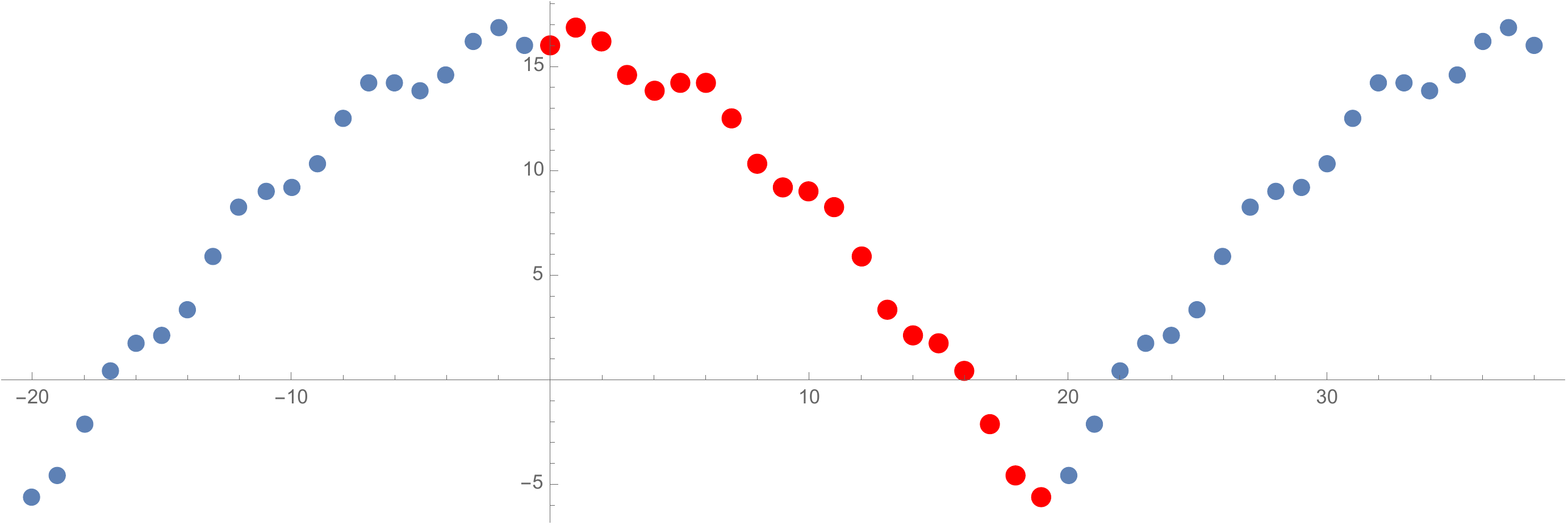}
\caption{(\textsc{top}) Periodic extension of data $s_0,s_1,\ldots,s_{19}$ (red).  The  extension belongs to $L^2(\Z/20\Z)$.
(\textsc{middle}) Periodic extension of the same data, but with odd boundary conditions.
This signal belongs to $L^2_-(\Z/40\Z)$.
(\textsc{bottom}) Periodic extension of the same data, but with even boundary 
conditions at both edges.  This signal belongs to $L^2_+(\Z/40\Z)$.
 Its ``smoothness'' suggests that the DCT may be of more practical use than the DST, or even the DFT,
 for the manipulation, storage, or compression of ``natural'' data.}
\label{Figure:Signals}
\end{figure}

Why is the DCT more prevalent than the DST?  
Suppose that we have samples
$s_0,s_1,\ldots,s_{m-1}$ taken at times $t=0,1,2\ldots,m-1$.
To employ Fourier-analytic techniques, this signal must be extended to $t \in \Z$ in a periodic fashion.  
For many applications, it behooves
the user to make this extension ``smooth'' in the sense that there are not large discrepancies 
between adjacent values.  This suggests the use of a reflection and even boundary conditions; see Figure \ref{Figure:Signals}.
A standard dictum in Fourier analysis is that greater smoothness of the input signal translates into
more rapid numerical convergence of associated algorithms.
The periodic extension of the sample that is used by the DCT is naturally ``smoother'' 
(for typical real-world signals) than those utilized by the DFT or DST.  Consequently, 
it is the DCT that plays a central role in modern signal processing.

There are many subtle variants of ``the'' DCT that appear in the literature, along with
their multidimensional analogues.  
Our particular selection is the most suitable from the viewpoint of supercharacter theory.
Indeed, our DCT matrix is precisely the $U$-matrix that
arises from a particularly simple supercharacter theory on $\Z/n\Z$.

Let $G = \Z/n\Z$ and $\zeta = \exp(-2\pi i /n)$.
Consider the action of the subgroup $\Gamma = \{ \pm 1\}$ of $(\Z/n\Z)^{\times}$ upon $G$. 
This produces the orbit decomposition
\begin{equation*}
		\X = 
		\begin{cases}
		\big\{ \{0\},\, \{\pm1\},\, \{ \pm 2\}, \ldots,\{\tfrac{n}{2}\pm 1\},\,  \{ \tfrac{n}{2} \} \big\} & \text{if $n$ is even},\\[5pt]
		\big\{ \{0\},\, \{\pm1\},\, \{ \pm 2\}, \ldots,\{\tfrac{n\pm 1}{2}\} \big\} & \text{if $n$ is odd}.
		\end{cases}
\end{equation*}
Let $N = |\X| = \floor{\frac{n}{2}}$. For $j = 0,1,\ldots,N$, we define
the corresponding superclasses
\begin{equation*}
X_j =
\begin{cases}
\{j,-j\} & \text{if } 2j \neq 0, \\
\{j\} & \text{if } 2j = 0.
\end{cases}
\end{equation*}
For $j=0,1,\ldots,N$, we have the supercharacters
\begin{equation}\label{sigma}
\sigma_j(k) = 
    \begin{cases}
        \zeta^{jk} + \zeta^{-jk} & \text{if $2j \neq 0$}, \\
        \zeta^{jk} & \text {if $2j=0$}.
        \end{cases}
\end{equation}
Euler's formula tells that 
\begin{equation}\label{eq:SigmaCos}
\sigma_j(\pm k) = |X_j| \cos \Big( \frac{2 \pi jk}{n} \Big ),
\end{equation}
in which $|X_j| \in \{1,2\}$ is the cardinality of $X_j$. 
\emph{We index the superclasses starting at $0$ rather than $1$. }
Doing so ensures that $i \in X_i$ for all $i \in G$, and so we may consider group elements and indices interchangeably. 
This convenience is more than enough to justify what is a small burden of notation. 

In the notation of Theorem \ref{Theorem:Big}, we have
\begin{equation}\label{UDCT}
[U]_{j+1,k+1}
=  \frac{\sqrt{|X_j||X_k|}}{\sqrt{n}} \cos \Big( \frac{2 \pi jk}{n} \Big),
\end{equation}
or more explicitly,
\begin{equation*}\small
U = \frac{1}{\sqrt{n}}\! \begin{bmatrix} 1 & \sqrt{2} & \sqrt{2} & \cdots & \sqrt{2} & 1 \\ \sqrt{2} & 2 \cos \frac{2 \pi}{n} &  2 \cos \frac{4 \pi}{n} & \cdots & 2 \cos \frac{(n-2) \pi}{n} & - \sqrt{2} \\ \sqrt{2} & 2 \cos \frac{4 \pi}{n} &  2 \cos \frac{8 \pi}{n} & \cdots & 2 \cos \frac{2(n-2) \pi}{n} & \sqrt{2} \\ \vdots & \vdots & \vdots & \ddots & \vdots & \vdots \\ \sqrt{2} & 2 \cos \frac{(n-2) \pi}{n} & 2 \cos \frac{2(n-2) \pi}{n} & \cdots & 2 \cos \frac{2(\frac{n}{2}-1)^2 \pi}{n} & (-1)^{\frac{n}{2}-1} \sqrt{2} \\ 1 & - \sqrt{2} & \sqrt{2} & \cdots & (-1)^{\frac{n}{2}-1} \sqrt{2} & (-1)^{\frac{n}{2}} 
\end{bmatrix}
\end{equation*}
if $n$ is even and
\begin{equation*}\small
 U = \frac{1}{\sqrt{n}}\! \begin{bmatrix} 1 & \sqrt{2}  & \sqrt{2} & \cdots & \sqrt{2} & \sqrt{2}  \\ \sqrt{2} & 2 \cos \frac{2 \pi}{n} &  2 \cos \frac{4 \pi}{n} & \cdots &  2 \cos \frac{(n-3) \pi}{n} &  2 \cos \frac{(n-1) \pi}{n} \\ \sqrt{2} & 2 \cos \frac{4 \pi}{n} &  2 \cos \frac{8 \pi}{n} & \cdots &  2 \cos \frac{2(n-3) \pi}{n} &  2 \cos \frac{2(n-1) \pi}{n} \\ \vdots & \vdots & \vdots & \ddots & \vdots & \vdots \\ \sqrt{2} & 2 \cos \frac{(n-3) \pi}{n} &  2 \cos \frac{2(n-3) \pi}{n} & \cdots &  2 \cos \frac{(n-3)^2 \pi}{n} &  2 \cos \frac{(n-3)(n-1) \pi}{n} \\ \sqrt{2} & 2 \cos \frac{(n-1) \pi}{n} &  2 \cos \frac{2(n-1) \pi}{n} & \cdots &  2 \cos \frac{(n-3)(n-1) \pi}{n} &  2 \cos \frac{(n-1)^2 \pi}{n} \end{bmatrix}
\end{equation*}
if $n$ is odd.  These so-called DCT matrices are real, symmetric, and unitary. 
They belong to $\M_{N+1}$, the set of $(N+1) \times (N+1)$ matrices.

The main result of this section identifies the matrices diagonalized by the DCT matrix \eqref{UDCT}.
Let $c_{ijk}$ denote the number of 
distinct solutions $(x,y) \in X_i \times X_j$ to 
$x + y = z$,
in which $z\in X_k$ is fixed.  As stated in Theorem \ref{Theorem:Big}, $c_{ijk}$
is independent of the particular representation $z\in X_k$ that is chosen.

\begin{theorem}\label{Theorem:DCT}
    Let $G = \Z/n\Z$, $N = \lfloor \frac{n}{2} \rfloor$, and let $U \in \M_{N+1}$ be the 
    discrete cosine transform matrix \eqref{UDCT}.
The matrices $T_0,T_1,\ldots,T_N \in \M_{N+1}$ defined by
\begin{equation}\label{eq:cijk}
[T_i]_{j+1,k+1} = \frac{c_{ijk} \sqrt{|X_k|}}{\sqrt{|X_j|}} 
\end{equation}
form a basis for the algebra $\A$ of matrices that are diagonalized by $U$.
They are real, symmetric, and satisfy
$$T_i = UD_iU^*,$$ in which
\begin{equation*}
    D_i = |X_i| \diag (  1,  \cos  \tfrac{2 \pi i}{n} , \cos \tfrac{4 \pi i}{n} , 
\ldots,  \cos \tfrac{2 \pi N i}{n}  ) \in \M_{N+1}.
\end{equation*}
 Moreover,  $T_0 = I$ and $T_i$ generates $\A$ if and only if $i$ is relatively prime to $n$.
       The most general matrix $T \in \M_{N+1}$ diagonalized by $U$ is 
            \begin{equation*}
            [T]_{j,k} = 
            \begin{cases}
                 t_{\min(n-|k-j|,|k-j|)} + t_{\min(n-k-j+2,k+j-2)} & \text{for $1 < j,k < \frac{n}{2} + 1$}, \\[10pt]
                 |X_{j-1}|^{\frac{1}{2}} |X_{k-1}|^{\frac{1}{2}} t_{\min(j-1,k-1)} & \text{for $j = 1$ or $k = 1$}, \\[10pt]
                 |X_{\tfrac{n}{2}+1-j}|^{\frac{1}{2}} |X_{\tfrac{n}{2}+1-k}|^{\frac{1}{2}}  t_{\max(\frac{n}{2}+1-j,\frac{n}{2}+1-k)} & 
                 \text{for $j = \frac{n}{2} +1$ or $k = \frac{n}{2} + 1$},
            \end{cases}
            \end{equation*}
        in which $t_0,t_1,\ldots,t_N\in\C$ are parameters
        (the last case only occurs if $n$ is even).
\end{theorem}

We defer the proof until Section \ref{Section:ProofDCT}.
Instead, we focus on several examples.

\begin{example}\label{Example:DCTeven}
    If $n$ is even, then $N = n/2$ and $T$ is
    \begin{equation*}\small
        \begin{bmatrix} 
            t_0 & \sqrt{2} t_1 & \sqrt{2}t_2 & \sqrt{2}t_3 & \cdots & \sqrt{2}t_{N-1} & t_N \\ 
            \sqrt{2} t_1 & t_0 + t_2 & t_1 + t_3 & t_2+t_4 & \cdots & t_{N-2} + t_N & \sqrt{2} t_{N-1} \\ 
            \sqrt{2} t_2 & t_{1} + t_3 & t_0 + t_4 & t_1+t_5 & \cdots & t_{N-3} + t_{N-1} & \sqrt{2} t_{N-2}  \\ 
            \sqrt{2} t_3 & t_2+t_4 & t_1+t_5 & t_0+t_6 & \cdots & t_{N-4} + t_{N-2} & \sqrt{2} t_{N-3} \\ 
            \vdots & \vdots & \vdots & \vdots & \ddots  & \vdots & \vdots \\ 
            \sqrt{2}t_{N-1} & t_{N-2} + t_N & t_{N-3} + t_{N-1} & t_{N-4}+t_{N-2} & \cdots & t_0 + t_2 & \sqrt{2}t_1 \\
            t_N & \sqrt{2}t_{N-1} & \sqrt{2}t_{N-2} & \sqrt{2}t_{N-3} & \cdots & \sqrt{2} t_{1} & t_0  
        \end{bmatrix} .
    \end{equation*}
\end{example}

\begin{example}\label{Example:DCTodd}
    If $n$ is odd, then $N = \floor{n/2}$ and $T$ is
    \begin{equation*}
        \small
        \begin{bmatrix} 
            t_0 & \sqrt{2} t_1 & \sqrt{2}t_2 & \sqrt{2} t_3 & \cdots & \sqrt{2}t_{N-1} & \sqrt{2} t_N \\ 
            \sqrt{2} t_1 & t_0 + t_2 & t_1 + t_3 & t_2+t_4 & \cdots & t_{N-2} + t_N & t_{N-1} + t_N \\ 
            \sqrt{2}t_2 & t_{1} + t_3 & t_0 + t_4 & t_1+t_5 & \cdots & t_{N-3} + t_{N-1} & t_{N-2} + t_{N-1}  \\ 
            \sqrt{2}t_3 & t_2+t_4 & t_1+t_5 & t_0 + t_6 & \cdots & t_{N-4}+t_{N-2} & t_{N-3} + t_{N-2} \\ 
            \vdots & \vdots & \vdots & \vdots & \ddots & \vdots & \vdots \\ 
            \sqrt{2}t_{N-1} & t_{N-2}+t_N & t_{N-3} + t_{N-1} & t_{N-4}+t_{N-2} & \cdots & t_0 + t_3 & t_1+t_2 \\
            \sqrt{2} t_N & t_{N-1}+t_N & t_{N-2} + t_{N-1} & t_{N-3} + t_{N-2} & \cdots & t_{1} + t_{2} & t_0 + t_{1}  
        \end{bmatrix}.
    \end{equation*}
\end{example}

The combinatorial aspect of Theorem \ref{Theorem:DCT} deserves special attention.

\begin{example}
    If $n=7$, then
    \begin{equation*}
        X_0 = \{0\}, \qquad X_1 = \{1,6\}, \qquad X_2 = \{2,5\},\quad \text{and} \quad X_3 = \{3,4\}.
    \end{equation*}
    The only solution in $X_3 \times X_1$ to $x + y = 3$ is $(4,6)$.
    Consequently, \eqref{eq:cijk} 
    produces
    \begin{equation*}
        [T_3]_{2,4} = \frac{c_{313} \sqrt{|X_3|}}{\sqrt{|X_1|}} = 1.
    \end{equation*}
    The two solutions in $X_3 \times X_3$ to $x+y=  0$ are
    $(3,4)$ and $(4,3)$.  Thus,
    \begin{equation*}
        [T_3]_{4,1} = \frac{c_{330} \sqrt{|X_0|}}{\sqrt{|X_3|}} = \sqrt{2} .
    \end{equation*}
    Computing the remaining entries in a similar fashion yields
    \begin{equation*}
        T_3 = \small
        \begin{bmatrix} 
            \0 & \0 & \0 & \sqrt{2} \\ 
            \0 & \0 & 1 & 1 \\ 
            \0 & 1 & 1 & \0 \\ 
            \sqrt{2} & 1 & \0 & \0 
        \end{bmatrix}.
    \end{equation*}
\end{example}

\begin{example}
    If $n=8$, then
    \begin{equation*}
        X_0 = \{0\}, \qquad X_1 = \{1,7\}, \qquad X_2 = \{2,6\}, \qquad X_3 = \{3,5\}, \quad \text{and} \quad X_4 = \{4\}.
    \end{equation*}
    The solutions in $X_3 \times X_1$ to $x+y =  4$ are $(3,1)$ and $(5,7)$.  Thus, \eqref{eq:cijk} produces
    \begin{equation*}
        [T_3]_{2,5} = \frac{c_{314} \sqrt{|X_4|}}{\sqrt{|X_1|}} = \sqrt{2}. 
    \end{equation*}
    The only solution in $X_3 \times X_4$ to $x + y = 1$ is $(5,4)$.  Thus,
    \begin{equation*}
        [T_3]_{5,2} = \frac{c_{341} \sqrt{|X_1|}}{\sqrt{|X_4|}} = \sqrt{2}. 
    \end{equation*}
    Computing the remaining entries in a similar fashion yields
    \begin{equation*}
        T_3 = \small
        \begin{bmatrix}
            \0 & \0 & \0 & \sqrt{2} & \0 \\
            \0 & \0 & 1 & \0 & \sqrt{2} \\
            \0 & 1 & \0 & 1 & \0 \\
            \sqrt{2} & \0 & 1 & \0 & \0 \\
            \0 & \sqrt{2} & \0 & \0 & \0
        \end{bmatrix}.
    \end{equation*}
\end{example}

\begin{example}
For $n = 10$, the most general matrix that is diagonalized by $U$ is
\begin{equation*}\small
    \begin{bmatrix} 
        t_0 & \sqrt{2} t_1 & \sqrt{2} t_2 & \sqrt{2} t_3 & \sqrt{2} t_4 & t_5 \\ 
        \sqrt{2} t_1 & t_0 + t_2 & t_1 + t_3 & t_2 + t_4 & t_3 + t_5 & \sqrt{2} t_4 \\
        \sqrt{2} t_2 & t_{1} + t_3 & t_0 + t_4 & t_1 + t_5 & t_2 + t_4 & \sqrt{2} t_3 \\ 
        \sqrt{2} t_3 & t_{2} + t_4 & t_{1} + t_5 & t_0 + t_4 & t_1 + t_3 & \sqrt{2} t_2 \\ 
        \sqrt{2} t_4 & t_{3} + t_5 & t_{2} + t_4 & t_{1} + t_3 & t_0 + t_2 & \sqrt{2} t_1 \\ 
        t_5 & \sqrt{2} t_4 & \sqrt{2} t_{3} & \sqrt{2} t_{2} & \sqrt{2} t_{1}  & t_0
    \end{bmatrix}
\end{equation*}
in which $t_0,t_1,t_2,t_3,t_4,t_5$ are free parameters.
It is a linear combination of 
\begin{align*}
    T_0 &= \small
    \begin{bmatrix} 
        1 & \0 & \0 & \0 & \0 & \0 \\ 
        \0 & 1 & \0 & \0 & \0 & \0 \\ 
        \0 & \0 & 1 & \0 & \0 & \0 \\ 
        \0 & \0 & \0 & 1 & \0 & \0 \\ 
        \0 & \0 & \0 & \0 & 1 & \0 \\ 
        \0 & \0 & \0 & \0 & \0 & 1 \\ 
    \end{bmatrix}
    &T_1 &= \small
    \begin{bmatrix} 
        \0 & \sqrt{2} & \0 & \0 & \0 & \0 \\ 
        \sqrt{2} & \0 & 1 & \0 & \0 & \0 \\ 
        \0 & 1 & \0 & 1 & \0 & \0 \\ 
        \0 & \0 & 1 & \0 & 1 & \0 \\ 
        \0 & \0 & \0 & 1 & \0 & \sqrt{2} \\ 
        \0 & \0 & \0 & \0 & \sqrt{2} & \0 
    \end{bmatrix},\\
    T_2&= \small
    \begin{bmatrix} 
        \0 & \0 & \sqrt{2} & \0 & \0 & \0 \\ 
        \0 & 1 & \0 & 1 & \0 & \0 \\ 
        \sqrt{2} & \0 & \0 & \0 & 1 & \0 \\ 
        \0 & 1 & \0 & \0 & \0 & \sqrt{2} \\ 
        \0 & \0 & 1 & \0 & 1 & \0 \\ 
        \0 & \0 & \0 & \sqrt{2} & \0 & \0 
    \end{bmatrix}, 
    &T_3 &= \small
    \begin{bmatrix} 
        \0 & \0 & \0 & \sqrt{2} & \0 & \0 \\ 
        \0 & \0 & 1 & \0 & 1 & \0 \\ 
        \0 & 1 & \0 & \0 & \0 & \sqrt{2} \\ 
        \sqrt{2} & \0 & \0 & \0 & 1 & \0 \\ 
        \0 & 1 & \0 & 1 & \0 & \0 \\ 
        \0 & \0 & \sqrt{2} & \0 & \0 & \0 
    \end{bmatrix},\\
    T_4 &= \small
    \begin{bmatrix} 
        \0 & \0 & \0 & \0 & \sqrt{2} & \0 \\ 
        \0 & \0 & \0 & 1 & \0 & \sqrt{2} \\ 
        \0 & \0 & 1 & \0 & 1 & \0 \\ 
        \0 & 1 & \0 & 1 & \0 & \0 \\ 
        \sqrt{2} & \0 & 1 & \0 & \0 & \0 \\ 
        \0 & \sqrt{2} & \0 & \0 & \0 & \0 
    \end{bmatrix},
    &T_5&= \small
    \begin{bmatrix} 
        \0 & \0 & \0 & \0 & \0 & 1 \\ 
        \0 & \0 & \0 & \0 & 1 & \0 \\ 
        \0 & \0 & \0 & 1 & \0 & \0 \\ 
        \0 & \0 & 1 & \0 & \0 & \0 \\ 
        \0 & 1 & \0 & \0 & \0 & \0 \\ 
        1 & \0 & \0 & \0 & \0  & \0
    \end{bmatrix} .
\end{align*}

\end{example}

\begin{example}
    For $n = 11$, the most general matrix that is diagonalized by $U$ is
    \begin{equation*}
    T = \small
    \begin{bmatrix} 
    t_0 & \sqrt{2} t_1 & \sqrt{2} t_2 & \sqrt{2} t_3 & \sqrt{2} t_4 & \sqrt{2} t_5 \\
    \sqrt{2} t_1 & t_0 + t_2 & t_1 + t_3 & t_2 + t_4 & t_3 + t_5 & t_4 + t_5\\ 
    \sqrt{2} t_2 & t_{1} + t_3 & t_0 + t_4 & t_1 + t_5 & t_2 + t_5 & t_3 + t_4 \\ 
    \sqrt{2} t_3 & t_{2} + t_4 & t_{1} + t_5 & t_0 + t_5 & t_1 + t_4 & t_2 + t_3\\ 
    \sqrt{2} t_4 & t_{3} + t_5 & t_{2} + t_5 & t_{1} + t_4 & t_0 + t_3 & t_1 + t_2 \\ 
    \sqrt{2} t_5 & t_{4} + t_5 & t_{3} + t_4 & t_{2} + t_3 & t_{1} + t_2 & t_0 + t_{1} 
    \end{bmatrix}
    \end{equation*}
    in which $t_0,t_1,t_2,t_3,t_4,t_5$ are free parameters.
    It is a linear combination of 
    \begin{align*}
        T_0 &= \small
        \begin{bmatrix} 
            1 & \0 & \0 & \0 & \0 & \0 \\ 
            \0 & 1 & \0 & \0 & \0 & \0 \\ 
            \0 & \0 & 1 & \0 & \0 & \0 \\ 
            \0 & \0 & \0 & 1 & \0 & \0 \\ 
            \0 & \0 & \0 & \0 & 1 & \0 \\ 
            \0 & \0 & \0 & \0 & \0 & 1 \\ 
        \end{bmatrix},
        &T_1&=\small
        \begin{bmatrix} 
            \0 & \sqrt{2} & \0 & \0 & \0 & \0 \\ 
            \sqrt{2} & \0 & 1 & \0 & \0 & \0 \\ 
            \0 & 1 & \0 & 1 & \0 & \0 \\ 
            \0 & \0 & 1 & \0 & 1 & \0 \\ 
            \0 & \0 & \0 & 1 & \0 & 1 \\ 
            \0 & \0 & \0 & \0 & 1 & 1 
        \end{bmatrix}, \\ 
        T_2&= \small
        \begin{bmatrix} 
            \0 & \0 & \sqrt{2} & \0 & \0 & \0 \\ 
            \0 & 1 & \0 & 1 & \0 & \0 \\ 
            \sqrt{2} & \0 & \0 & \0 & 1 & \0 \\ 
            \0 & 1 & \0 & \0 & \0 & 1 \\ 
            \0 & \0 & 1 & \0 & \0 & 1 \\ 
            \0 & \0 & \0 & 1 & 1 & \0 
        \end{bmatrix},
        &T_3&= \small
        \begin{bmatrix} 
            \0 & \0 & \0 & \sqrt{2} & \0 & \0 \\ 
            \0 & \0 & 1 & \0 & 1 & \0 \\ 
            \0 & 1 & \0 & \0 & \0 & 1 \\ 
            \sqrt{2} & \0 & \0 & \0 & \0 & 1 \\ 
            \0 & 1 & \0 & \0 & 1 & \0 \\ 
            \0 & \0 & 1 & 1 & \0 & \0 
        \end{bmatrix}, \\
        T_4 &= \small
        \begin{bmatrix} 
        \0 & \0 & \0 & \0 & \sqrt{2} & \0 \\ 
        \0 & \0 & \0 & 1 & \0 & 1 \\ 
        \0 & \0 & 1 & \0 & \0 & 1 \\ 
        \0 & 1 & \0 & \0 & 1 & \0 \\  
        \sqrt{2} & \0 & \0 & 1 & \0 & \0 \\ 
        \0 & 1 & 1 & \0 & \0 & \0 
        \end{bmatrix},
        &T_5 &= \small
        \begin{bmatrix} 
            \0 & \0 & \0 & \0 & \0 & \sqrt{2} \\ 
            \0 & \0 & \0 & \0 & 1 & 1 \\ 
            \0 & \0 & \0 & 1 & 1 & \0 \\ 
            \0 & \0 & 1 & 1 & \0 & \0 \\ 
            \0 & 1 & 1 & \0 & \0 & \0 \\ 
            \sqrt{2} & 1 & \0 & \0 & \0 & \0 
        \end{bmatrix}.
    \end{align*}
    The matrices above are analogous to those encountered by Feig and Ben-Or \cite{ADCT}, 
    who considered the modified DCT matrix 
    \begin{equation*}
        [C_n]_{i,j} = c_i \cos \frac{ 2 \pi (2j-1)(i-1)}{4n} ,
    \end{equation*}
    in which $c_i = \sqrt{1/n}$ for $i=1$ and $\sqrt{2/n}$ otherwise. 
\end{example}

\begin{example}
Matrices diagonalized by the DCT have been studied before, 
but with different techniques and sometimes with different DCT matrices \cite{DPDCT,ADCT}.  
Theorem \ref{Theorem:DCT} recovers many established results.
For example, the matrix 
\begin{equation*}\small
    \begin{bmatrix} 
        \0 & 1/\sqrt{2} & \0 & \cdots & \cdots & \0 \\ 
        1/\sqrt{2} & \ddots & 1/2 & \ddots & & \vdots \\
        0 & 1/2 & \ddots & \ddots & \ddots & \vdots \\ 
        \vdots & \ddots & \ddots & \ddots & 1/2 & \0 \\ 
        \vdots & & \ddots & 1/2 & \ddots & 1/\sqrt{2} \\ 
        \0 & \cdots & \cdots & \0 & 1/\sqrt{2} & \0 
    \end{bmatrix},
\end{equation*}
appears in \cite{DPDCT}.  
In our notation, it corresponds to even $n$ and parameters
$t_0 = 0$, $t_1 = 1/2$, and $t_2 =t_3= \cdots = t_N = 0$; see Example \ref{Example:DCTeven}.
\end{example}

\begin{example}
For $n$ odd, the bottom right $N \times N$ submatrix of any matrix diagonalized by $U$ is 
a Toeplitz plus Hankel matrix:
\begin{equation*}
\begin{bmatrix} 
t_0 & t_1 & \cdots & t_{N-1} \\ 
t_1 & \ddots & \ddots & \vdots \\ 
\vdots & \ddots & \ddots & t_1 \\ 
t_{N-1} & \cdots & t_1 & t_0 
\end{bmatrix} 
+ 
\begin{bmatrix} 
t_{2} & \cdots & t_{N} & t_{N} \\ 
\vdots & \iddots & \iddots & t_{N-1} \\ 
t_{N} & \iddots & \iddots & \vdots \\ 
t_N & t_{N-1} & \cdots & t_1 
\end{bmatrix} 
\end{equation*}
An analogous presentation exists when $n$ is even if we also exclude the first and last row and column. 
In \cite{DPDCT} it is shown that the DCT-I matrix, obtained by replacing $|X_k|^{1/2}$ and $|X_j|^{1/2}$ 
with $|X_k|$ and $|X_j|$, respectively, in (\ref{UDCT}), diagonalizes matrices that are genuinely Toeplitz plus Hankel. In \cite{FSD} Grishin and Strohmer demonstrate that it is simple to go from the DCT-I to $U$, and that there are advantages to both matrices. 
While the DCT-I diagonalizes certain Toeplitz plus Hankel matrices, it is not unitary like $U$.
\end{example}

Theorem \ref{Theorem:DCT} also provides an explanation for this Toeplitz plus Hankel structure. 
The matrix entry $[T_i]_{j+1,k+1} = c_{ijk}$ is nonzero if and only if $\pm i \pm j = k$, or equivalently, when $i \in X_{k-j}$ or $i \in X_{k+j}$. 
If $i \in X_{k-j}$, then $i \in X_{(k+\ell)-(j+\ell)}$ for any $\ell \in G$. 
Thus, along the diagonal that contains $(j+1,k+1)$, 
$T_i$ is always nonzero; this gives us one sub- or super-diagonal of a Toeplitz matrix. 
Similarly, if $i \in X_{k-j}$, then $T_i$ is nonzero along the entire anti-diagonal containing $(j,k)$, 
giving us a component of a Hankel matrix.  See \cite{Zellini} for a displacement-rank approach to such matrices.

\section{Proof of Theorem \ref{Theorem:DCT}}\label{Section:ProofDCT}

Let $\A$ denote the commutative, complex algebra of matrices that are diagonalized by $U$.
The algebra of $(N+1) \times (N+1)$ diagonal matrices has dimension $N+1$.  Thus, $\dim \A = N+1$.
The diagonal matrices $D_0,D_1,D_2,\ldots,D_N$ are linearly independent because their diagonals
$$[ 1\,\,  \cos  \tfrac{2 \pi i}{n} \,\, \cos \tfrac{4 \pi i}{n} \,\, \ldots \,\,  \cos \tfrac{2 \pi N i}{n}]^{\T} \in \C^{N+1}$$
are the columns of the matrix $[\sigma_{i-1}(j-1)]_{i,j=1}^{N+1}$, which is similar to the unitary matrix $U$.
Thus, $\{T_0,T_1,\ldots,T_N\}$ is linearly independent and hence it spans $\A$.

The eigenvalues $$1,\quad  \cos  \tfrac{2 \pi i}{n} ,\quad \cos \tfrac{4 \pi i}{n} , \ldots, \quad \cos \tfrac{2 \pi N i}{n}$$
of $T_i$ are distinct if and only if $i$ is relatively prime to $n$.  In this case,
the Lagrange interpolation theorem ensures that for any diagonal matrix $D \in \M_{N+1}$, there is a polynomial $p$
so that $p(T_i) = UDU^*$.  Thus, $T_i$ generates $\A$.

We claim that $T_0 = I$.
If $i = 0$, then $X_i = \{0\}$.  Consequently, $x+y \in X_k$ and $(x,y) \in X_i \times X_j$ and imply $j=k$;
moreover, $c_{ijj} = 1$.  Thus, $T_0 = I$.

We now consider $T_i$ for $i=1,2,\ldots,N$ and identify the locations of all nonzero entries in each matrix.
First suppose that $n$ is odd (if $n$ is even then there are a few additional cases to consider; we will do this later). 

If $j = 0$, then the argument above implies that $i = k$. 
Since $n$ is odd, $-i = k$ means $i = 0$. 
Thus, $c_{ijk} = 1$ and, since $k = i \neq 0$, we have $|X_k| = 2$. 
By symmetry,
\begin{equation*}
     [T_i]_{i+1,1} = [T_i]_{1,i+1} = \sqrt{2}.
\end{equation*}
An analogous approach applies if $k = 0$. In all other cases, 
$i,j,k$ are nonzero and hence $|X_j|=|X_k|=2$. 
\begin{enumerate}
    \item[(i)] Suppose that $c_{ijk} = 2$.  Without loss of generality, let $(i,j)$ be one of the solutions
    to $x+y = k$ with $(x,y) \in X_i\times X_j$.  The other potential solution must be
    one of $(i,-j)$, $(-i,-j)$, or $(-i,j)$.  These possibilities imply that 
    $2j = 0$, $2k=0$, or $2i = 0$, respectively. Since $i,j,k\neq 0$, this is not possible.
    \item[(ii)] 
    Suppose that $c_{ijk} = 1$, with $(i,j)$ as the solution.  
    We see that $\pm i \pm j = k$ if and only if $i \in X_{j+k}$ or $i \in X_{k-j}$. For such $i$,
    \begin{equation*}
    [T_i]_{j+1,k+1} = 1.
    \end{equation*}
    Since $c_{ijk} \in \{0,1,2\}$, it follows that $T_i$ is $0$ elsewhere. 
\end{enumerate}

If $T \in \A$, then $T = \sum_{i=0}^N t_i T_i$ for some $t_0,t_1,\ldots,t_N \in \C$.
The preceding analysis implies that $T$ equals
\begin{equation*} \footnotesize
    \begin{bmatrix} 
        t_0 & \sqrt{2} t_1 & \sqrt{2}t_2 & \sqrt{2} t_3 & \cdots & \sqrt{2} t_{N-1} & \sqrt{2} t_N \\ 
        \sqrt{2} t_1 & t_0 + t_2 & t_1 + t_3 & t_2+t_4 & \cdots & t_{N-2} + t_N & t_{N-1} + t_{N+1} \\ 
        \sqrt{2}t_2 & t_{-1} + t_3 & t_0 + t_4 & t_1+t_5 & \cdots & t_{N-3} + t_{N+1} & t_{N-2} + t_{N+2}  \\ 
        \sqrt{2}t_3 & t_{-2}+t_4 & t_{-1}+t_5 & t_0 + t_6 & \cdots & t_{N-4}+t_{N+3} & t_{N-3} + t_{N+3} \\ 
        \vdots & \vdots & \vdots & \vdots & \ddots & \vdots & \vdots \\ 
        \sqrt{2}t_{N-1} & t_{2-N}+t_N & t_{3-N} + t_{N+1} & t_{4-N}+t_{N+2} & \cdots & t_0 + t_{2N-2} & \sqrt{2} t_1 \\
        \sqrt{2} t_N & t_{1-N}+t_{N+1} & t_{2-N}+t_{N+2} & t_{3-N}+t_{N+3} & \cdots & t_{-1} + t_{2N-1} & t_0 + t_{2N}  
    \end{bmatrix} 
\end{equation*}
in which, for the sake of convenience, we let
$t_i = t_{-i} = t_{n-i}$ for all $i$.
The preceding simplifies to the matrix presented in Example \ref{Example:DCTodd}.

Now suppose that $n$ is even.  The preceding results largely carry over, but there are now extra cases to consider.
\begin{enumerate}
    \item[(iii)] Suppose that $i = N = \frac{n}{2}$. Then the only solutions to $x + y \in X_k$ with $(x,y) \in X_i\times X_j$
     are when, without loss of generality, $k = N-j$.
     Since $|X_j| = |X_{N-j}|$ for all $j$, an appeal to \eqref{eq:cijk} reveals that $T_N$ is the reversed identity matrix.

    \item[(iv)] Suppose that $i \neq 0$, $i \neq N$, and $c_{ijk} = 2$.
    In addition to the cases identified in (ii), 
    we now have the possibilities $j = N$ or $k = N$. From \eqref{eq:cijk} we obtain
    \begin{equation*}
        [T_i]_{N+1,N-i+1} = [T_i]_{N-i+1,N+1} = \sqrt{2}.
    \end{equation*}
\end{enumerate}
In all other cases, $i,j,k \not \in \{0,N\}$, so the rest of our analysis from the odd case
carries over.  
If $T \in \A$, then $T = \sum_{i=0}^N t_i T_i$ for some $t_0,t_1,\ldots,t_N \in \C$.
The preceding analysis implies that $T$ equals
\begin{equation*}\footnotesize
        \begin{bmatrix} 
            t_0 & \sqrt{2} t_1 & \sqrt{2}t_2 & \sqrt{2}t_3 & \cdots & \sqrt{2}t_{N-1} & t_N \\ 
            \sqrt{2} t_1 & t_0 + t_2 & t_1 + t_3 & t_2+t_4 & \cdots & t_{N-2} + t_N & \sqrt{2} t_{N-1} \\ 
            \sqrt{2} t_2 & t_{-1} + t_3 & t_0 + t_4 & t_1+t_5 & \cdots & t_{N-3} + t_{N-1} & \sqrt{2} t_{N+2}  \\ 
            \sqrt{2} t_3 & t_{-2}+t_4 & t_{-1}+t_5 & t_0+t_6 & \cdots & t_{N-4} + t_{N-2} & \sqrt{2} t_{N+3} \\ 
            \vdots & \vdots & \vdots & \vdots & \ddots  & \vdots & \vdots \\ 
            \sqrt{2}t_{N-1} & t_{2-N} + t_N & t_{3-N} + t_{N-1} & t_{4-N}+t_{N-2} & \cdots & t_0 + t_{2N-2} & \sqrt{2}t_1 \\
            t_N & \sqrt{2}t_{N-1} & \sqrt{2}t_{N-2} & \sqrt{2}t_{N-3} & \cdots & \sqrt{2} t_{1} & t_0  
        \end{bmatrix} .
\end{equation*}
This simplifies to the matrix presented in Example \ref{Example:DCTeven}. 

We can now obtain an explicit formula for the entries of the most general matrix $T = [T_{j,k}] \in \M_{N+1}$
diagonalized by $U$. The first row and column (and the last row and column if $n$ is even) have a different structure from the rest of the matrix; one can see that 
\begin{equation*} 
[T]_{j,k} =
\begin{cases} 
|X_{j-1}|^{\frac{1}{2}} |X_{k-1}|^{\frac{1}{2}} t_{\min(j-1,k-1)} & \text{for $j = 1$ or $k = 1$}, \\[8pt]
|X_{\tfrac{n}{2}+1-j}|^{\frac{1}{2}} |X_{\tfrac{n}{2}+1-k}|^{\frac{1}{2}}  t_{\max(\frac{n}{2}+1-j,\frac{n}{2}+1-k)} 
& \text{for $j = \frac{n}{2} +1$ or $k = \frac{n}{2} + 1$},
\end{cases} 
\end{equation*}
holds for these sections of $T$; the second case occurs only if $n$ is even.

We direct our attention now to the remaining entries. 
First observe that
    \begin{equation*}
        [T]_{j,k} = t_{k-1-(j-1)} + t_{k-1+j-1} = t_{k-j} + t_{k+j-2} 
    \end{equation*}
To ensure that all of our subscripts are between $0$ and $n-1$, we take the absolute value of the first subscript.
Since $1 < j,k < N < \frac{n}{2}$, it follows that $0 \leq k+j-2 \leq n-1$.  Thus,
    \begin{equation*}
    [T]_{j,k} = t_{|k-j|} + t_{k+j-2}.
    \end{equation*}
Finally, we need to ensure that our subscripts are between $0$ and $N$. 
If $N+1 \leq i\leq n-1$, then $n-i$ gives the correct index and is in the desired range. 
Consequently,
\begin{equation*}
[T]_{j,k} = t_{\min(n-|k-j|,|k-j|)} + t_{\min(n-k-j+2,k+j-2)} . \qed
\end{equation*}

\section{Discrete sine transform}
The discrete sine transform (DST) is the oft-neglected sibling of the DCT.
Since $L^2(G)$ is finite dimensional and $L^2_-(G)  =  L^2_+(G) ^{\perp}$, it follows from 
the DFT-invariance of $ L^2_+(G) $ that $L^2_-(G) $ is also DFT-invariant (recall that the DFT is a unitary operator).
As mentioned in Section \ref{Section:DCT}, the DST is the restriction of
the DFT to $L^2_-(G)$, the subspace of odd functions in $L^2(G)$.
Here $G = \Z/n\Z$, as usual.
Let
\begin{equation*}
N = \Big\lfloor\frac{n-\frac{1}{4}}{2} \Big\rfloor.
\end{equation*}
As before, the sets $X_j = \{j,-j\}$ for $j = 1,2,\ldots,N$, along with $\{0\}$ (and $\{n/2\}$ if $n$ is even), partition $G$. 

In our consideration of the DCT, we saw that the supercharacters \eqref{eq:SigmaCos}
are constant on each superclass.  In contrast, the corresponding ``supercharacters''
obtained by replacing cosines with sines are no longer constant on each superclass.
This is a crucial distinction between the DCT and DST:
the DST does not arise directly from  a supercharacter theory on $G$.
Nevertheless, we are still able to obtain an analogue of Theorem \ref{Theorem:DCT}
for the DST by appealing to the DFT-invariance of $L^2_-(G)$ and considering
the ``orthogonal complement'' of the DCT supercharacter theory.

Define $\tau_j(k) = \zeta^{-jk} - \zeta^{jk}$ for $j = 0,1,\ldots,N$.
Then $\{\tau_j\}_{j=1}^N$ is an orthogonal basis for $L^2_-(G) $ and
\begin{equation}\label{eq:Tau}
\tau_j(k) = -2i \sin \Big( \frac{2 \pi jk}{n} \Big), 
\end{equation}
in which $i$ denotes the imaginary unit. 
Normalizing the $\tau_j$ yields
\begin{equation*}
v_j(k) = \frac{\tau_j(k)}{\sqrt{2n}} = \frac{\sqrt{2} \sin(\frac{2 \pi jk}{n})}{i \sqrt{n}}.
\end{equation*}
Let $V_n \in \M_N$ denote the matrix representation of the restriction of $\F$ to $L^2_-(G) $ with respect to 
the orthonormal basis $\{v_j\}_{j=1}^N$.   Then $V_n$ is unitary and a computation confirms that
\begin{equation}
[V_n]_{j,k} 
= \frac{2}{i \sqrt{n}} \sin \Big( \frac{2 \pi jk}{n} \Big).  \label{UDST} 
\end{equation} 
Thus, 
\begin{equation*}
    V_n 
    =  \frac{2}{i \sqrt{n}} \small 
    \begin{bmatrix}  
    \sin \frac{2 \pi}{n} &  \sin \frac{4 \pi}{n} & \cdots &  \sin \frac{2N \pi}{n} \\[3pt]  
    \sin \frac{4 \pi}{n} &  \sin \frac{8 \pi}{n} & \cdots &  \sin \frac{4N \pi}{n} \\
    \vdots & \vdots & \ddots & \vdots \\[3pt]  
    \sin \frac{2N \pi}{n} &  \sin \frac{4N \pi}{n} & \cdots &  \sin \frac{2N^2 \pi}{n} 
    \end{bmatrix}
    .
\end{equation*} 
The matrices $V_n$ are purely imaginary, complex symmetric, and unitary. 
If $n$ is clear from context, we often omit the subscript and write $V$. 
\emph{Although the DST cannot be attacked directly via supercharacter theory}, we can use the DFT invariance of $ L^2_-(G) $
to obtain a satisfying analogue of Theorem \ref{Theorem:DCT}.

\begin{theorem}\label{Theorem:DST}
Let $G = \Z/n\Z$, $N = \floor{\frac{n-\frac{1}{4}}{2}}$, and let $V \in \M_N$ be the discrete sine transform matrix \eqref{UDST}. 
Let $\sgn x$ denote the sign of $x$; let $\sgn 0 = 0$.
\begin{enumerate}
\item The most general $S \in \M_N$ diagonalized by $V$ is given by 
\begin{equation}\label{eq:genDST} 
[S]_{j,k} = \sum_{\ell=1}^{\min(j,k)} \sgn\big(\tfrac{n}{2}-|k-j|-2\ell+1\big) s_{\min(n-|k-j|-2\ell+1,|k-j|+2\ell-1)},
\end{equation} 
in which $s_0  = 0$, and $s_1,s_2,\ldots,s_N\in \C$ 
are free parameters that correspond, in that order, to the entries in the first row of $S$. 
For $i=1,2,\ldots,N$, the matrices $S_i$ obtained by
setting $s_j = \delta_{i,j}$ in \eqref{eq:genDST}
form a basis for the algebra $\A$ diagonalized by $V$. In particular, $S_1 = I$.

\item Let $x_{i,j} = 1$ if $j \in X_i$ and 0 otherwise.
The matrices $T_1,T_2,\ldots,T_N \in \M_{N}$ defined by
\begin{equation}\label{eq:xijk} 
[T_i]_{j,k} = x_{i,j-k} - x_{i,j+k} 
\end{equation} 
are real, symmetric, and satisfy
$$T_i = VD_iV^*,$$
in which
\begin{equation*}
D_i = 2 \diag \big(\cos \tfrac{2 \pi i}{n}, \cos \tfrac{4 \pi i}{n}, \ldots, \cos \tfrac{2 \pi N i}{n} \big) \in \M_N.
\end{equation*}
Moreover,
$T_i$ generates $\A$ if and only $i$ is relative prime to $n$. 

\item If $n$ is odd, then $\{T_1,T_2,\ldots,T_N\}$ is a basis for $\A$.  
Another formula for the entries for a general $T \in \A$ is given by
\begin{equation}\label{eq:oddDST}
[T_i]_{j,k} = t_{\min(n-j-k,j+k)} - t_{\min(n-j+k,j-k)},
\end{equation}
in which $t_0 = 0$, and $t_1,t_2,\ldots,t_N\in \C$ are free parameters.
\end{enumerate}
\end{theorem}

For odd $n$, Theorem \ref{Theorem:DST} provides two bases for $\A$. 
The basis described in (a) is obtained by a brute force method which, if applied to the DCT, yields the basis in Theorem \ref{Theorem:DCT}. However, it is cumbersome to work with; the following examples illustrate its inelegance and unwieldiness. 
The basis obtained in (b) is superior in several ways.  Not only is it much simpler in appearance, it also has
a nice combinatorial explanation.

The matrices given by \eqref{eq:genDST} and \eqref{eq:oddDST} are easier to grasp with examples. 
We defer the proof of Theorem \ref{Theorem:DST} until Section \ref{Section:ProofDST} and focus on some instructive examples. 

\begin{example}\label{ex:DSTeven}
If $n$ is even, the most general matrix diagonalized by $V_n$ is
\begin{equation*}\small
    \begin{bmatrix}
    s_1 & s_2 & s_3 & \cdots & s_{N-1} & s_N \\ 
    s_2 & s_1 + s_3 & s_2 + s_4 & \cdots & s_{N-2}+s_N & s_{N-1} \\ 
    s_3 & s_2 + s_4 & s_1 + s_3 + s_5 & \cdots & s_{N-3}+s_{N-1} & s_{N-2} \\ 
    \vdots & \vdots & \vdots & \ddots & \vdots & \vdots \\ 
    s_{N-1} & s_{N-2}+s_N & s_{N-3}+s_{N-1} & \cdots & s_1+s_3 & s_2 \\ 
    s_N & s_{N-1} & s_{N-2} & \cdots & s_2 & s_1
    \end{bmatrix}
\end{equation*}
in which $s_1,s_2,\ldots,s_N$ are free parameters. 
Bini and Capovani were the first to call the matrix above a $\Tau$-class matrix, and referred to the algebra $\A$ as $\Tau_N$. 
This class of matrices occurs in the study of Toeplitz matrices and is known to be diagonalized by our DST matrix \cite{SCP}. We recapture this result, and with our method we are able to find an analogous basis for the case where $n$ is odd, which has been much less studied.  
These matrices also form a subspace of the Toeplitz plus Hankel matrices \cite{FTTE}. 
There is a considerable amount of literature on $\Tau$-class matrices because of
their desirable computational properties. For instance, a $\Tau_N$ matrix system can be solved in $O(N \log N)$ 
time using algorithms for centrosymmetric Toeplitz plus Hankel matrices \cite{CHM}.
This makes $\Tau$-class matrices suitable as preconditioners for banded Toeplitz systems \cite{FTTE, STBP, FDSM}.

From \cite{FDSM}, $\Tau$-class matrices may also be defined as the $N \times N$ matrices $A = [a_{ij}]_{i,j=1}^N$ whose entries satisfy the ``cross-sum" condition 
\begin{equation}\label{eq:xsum}
    a_{i-1,j} + a_{i+1,j} = a_{i,j-1} + a_{i,j+1} ,
\end{equation}
in which $a_{N+1,j} = a_{i,N+1} = a_{0,j} = a_{i,0} = 0$. 
\end{example}

\begin{example}\label{ex:DSTodd}
If $n$ is odd, the most general matrix that is diagonalized by $V_n$ is
\begin{equation*}
    \small \begin{bmatrix}
    s_{1} & s_{2} & s_3 & \cdots & s_{N-1} & s_{N} \\ 
    s_{2} & s_{1} + s_{3} & s_2+s_4 & \cdots & s_{N-2}+s_N & s_{N-1}+s_{N+1} \\ 
    \vdots & \vdots & \ddots & \vdots & \vdots & \vdots  \\ 
    \vdots & \vdots & \vdots & \ddots &  \vdots & \vdots \\ 
    s_{N-1} & s_{N-2}+s_N & s_{N-3}+s_{N-1}+s_{N+1} & \vdots &  \sum_{j=1}^{N-1} s_{2j-1} & \sum_{j=1}^{N-1} s_{2j} \\ 
    s_{N} & s_{N-1}+s_{N+1} + s_{N+1} & s_{N-2}+s_N+s_{N+2} & \cdots & \sum_{j=1}^{N-1} s_{2j} & \sum_{j=1}^N s_{2j-1}
    \end{bmatrix}
\end{equation*}
in which $s_1,s_2,\ldots,s_N$ are free parameters, 
and $s_i = -s_{n-i}$. 
A glance at Example \ref{ex:DSTeven} confirms that the even and odd cases are
strikingly different. Because of this unexpected complexity, the odd case, as mentioned in the preceding example, does not appear to have been addressed
completely in the literature before.

However, these matrices enjoy many of the same properties $\Tau$ matrices do; they are Toeplitz plus Hankel, symmetric, and diagonalized by the DST matrix \eqref{UDST}. Further, the same equation (\ref{eq:genDST}) used to obtain these matrices recovers the $\Tau$ matrices if $n$ is even, 
so we may consider \eqref{eq:genDST} as providing a generalization of $\Tau$ matrices.
Using \eqref{eq:oddDST}, a more transparent description is
\begin{equation*}\small
    T = \begin{bmatrix}
    t_2 & t_3 - t_1 & t_4 - t_2 & \cdots & t_N - t_{N-2} & t_N - t_{N-1} \\
    t_3 - t_1 & t_4 & t_5 - t_1 & \cdots & t_N - t_{N-3} & t_{N-1} - t_{N-2} \\
    t_4 - t_2 & t_5 - t_1 & t_6 & \cdots & t_{N-1} - t_{N-4} & t_{N-2} - t_{N-3} \\
    \vdots & \vdots & \vdots & \ddots & \vdots & \vdots \\
    t_N - t_{N-2} & t_N - t_{N-3} & t_{N-1}-t_{N-4} & \cdots & t_3 & t_2 - t_1 \\
    t_N - t_{N-1} & t_{N-1} - t_{N-2} & t_{N-2} - t_{N-3} & \cdots & t_2 - t_1 & t_1 
    \end{bmatrix}
\end{equation*}
in which $t_1,t_2,\ldots,t_N$ are free parameters. From this parameterization we see these matrices even almost satisfy 
\eqref{eq:xsum}, failing to hold only at the right edge.  For instance, considering the $(2,N)$ entry, 
\begin{equation*}
    [T]_{1,N} + T_{3,N} = t_N - t_{N-1} + t_{N-2} - t_{N-3} \neq [T]_{2,N-1} + [T]_{2,N+1} = t_N - t_{N-3}
\end{equation*}
since the cross-sum condition takes $[T]_{2,N+1} = 0$. 
\end{example}

\begin{example}
For $n = 11$, the most general matrix diagonalized by $V_n$ is
\begin{equation}\label{eq:S11VnTH}  \small
    \begin{bmatrix} 
    s_1 & s_2 & s_3 & s_4 & s_5 \\ 
    s_2 & s_1 + s_3 & s_2 + s_4 & s_3 + s_5 &  s_4 - s_{5} \\ 
    s_3 & s_2 + s_4 & s_1 + s_3 + s_5 & s_2 + s_4 - s_{5} & s_3 + s_5 - s_{4} \\ 
    s_4 & s_3 + s_5 & s_2 + s_4 - s_{5} & s_1 + s_3 + s_5 - s_{4} & s_2 + s_4 - s_5 - s_3 \\ 
    s_5 & s_4 - s_5 & s_3 + s_5 - s_4 & s_2 + s_4 - s_5 - s_3 & s_1 + s_3 + s_5 - s_4 - s_2 
    \end{bmatrix}
\end{equation}
in which $s_1,s_2,s_3,s_4,s_5\in \C$ are free parameters.
It is a linear combination of
\begin{equation*}\small
    S_1 = 
    \begin{bmatrix} 
        1 & \0 & \0 & \0 & \0 \\ 
        \0 & 1 & \0 & \0 & \0 \\ 
        \0 & \0 & 1 & \0 & \0 \\ 
        \0 & \0 & \0 & 1 & \0 \\ 
        \0 & \0 & \0 & \0 & 1 \\ 
    \end{bmatrix},
    \quad
    S_2 = 
    \begin{bmatrix} 
        \0 & 1 & \0 & \0 & \0 \\
        1 & \0 & 1 & \0 & \0 \\
        \0 & 1 & \0 & 1 & \0 \\
        \0 & \0 & 1 & \0 & 1 \\
        \0 & \0 & \0 & 1 & -1 
    \end{bmatrix},
    \quad
    S_3 = 
    \begin{bmatrix}
        \0 & \0 & 1 & \0 & \0 \\
        \0 & 1 & \0 & 1 & \0 \\
        1 & \0 & 1 & \0 & 1 \\
        \0 & 1 & \0 & 1 & -1 \\
        \0 & \0 & 1 & -1 & 1 
    \end{bmatrix}, \\
\end{equation*}
\begin{equation*}\small
    S_4  = 
    \begin{bmatrix}
    \0 & \0 & \0 & 1 & \0 \\
    \0 & \0 & 1 & \0 & 1 \\
    \0 & 1 & \0 & 1 & -1 \\
    1 & \0 & 1 & -1 & 1 \\
    \0 & 1 & -1 & 1 & -1
    \end{bmatrix},
    \quad\text{and}\quad
    S_5 =
    \begin{bmatrix}
    \0 & \0 & \0 & \0 & 1 \\
    \0 & \0 & \0 & 1 & -1 \\
    \0 & \0 & 1 & -1 & 1 \\
    \0 & 1 & -1 & 1 & -1 \\
    1 & -1 & 1 & -1 & 1 
    \end{bmatrix}.
\end{equation*}
It is apparent each $S_i$ is Toeplitz plus Hankel; hence \eqref{eq:S11VnTH} is Toeplitz plus Hankel as well. 
Using \eqref{eq:oddDST} we obtain the alternate parametrization
\begin{equation*}\small
    \begin{bmatrix}
    t_2 & t_3 - t_1 & t_4 - t_2 & t_5 - t_3 & t_5 - t_4 \\
    t_3 - t_1 & t_4 & t_5 - t_1 & t_5 - t_2 & t_4 - t_3 \\
    t_4 - t_2 & t_5 - t_1 & t_5 & t_4 - t_1 & t_3 - t_2 \\
    t_5 - t_3 & t_5 - t_2 & t_4 - t_1 & t_3 & t_2 - t_1 \\
    t_5 - t_4 & t_4 - t_3 & t_3 - t_2 & t_2 - t_1 & t_1 
    \end{bmatrix}
\end{equation*}
in which $t_1,t_2,t_3,t_4,t_5\in \C$ are free parameters.
It is a linear combination of
\begin{equation*}
        T_1 = \small
        \begin{bmatrix}
        \0 & -1 & \0 & \0 & \0 \\
        -1 & \0 & -1 & \0 & \0 \\
        \0 & -1 & \0 & -1 & \0 \\
        \0 & \0 & -1 & \0 & -1 \\
        \0 & \0 & \0 & -1 & 1 \\
        \end{bmatrix},\qquad\normalsize
        T_2= \small
        \begin{bmatrix}
        1 & \0 & -1 & \0 & \0 \\
        \0 & \0 & \0 & -1 & \0 \\
        -1 & \0 & \0 & \0 & -1 \\
        \0 & -1 & \0 & \0 & 1 \\
        \0 & \0 & -1 & 1 & \0 
        \end{bmatrix},
\end{equation*}
\begin{equation*}
        T_3  = \small
        \begin{bmatrix}
        \0 & 1 & \0 & -1 & \0 \\
        1 & \0 & \0 & \0 & -1 \\
        \0 & \0 & \0 & \0 & 1 \\
        -1 & \0 & \0 & 1 & \0 \\
        \0 & -1 & 1 & \0 & \0 
        \end{bmatrix}\!\!,\, \normalsize
        T_4  = \small
        \begin{bmatrix}
        \0 & \0 & 1 & \0 & -1 \\
        \0 & 1 & \0 & \0 & 1 \\
        1 & \0 & \0 & 1 & \0 \\
        \0 & \0 & 1 & \0 & \0 \\
        -1 & 1 & \0 & \0 & \0
        \end{bmatrix}\!\!,\, 
        \normalsize
            T_5 = \small
        \begin{bmatrix}
        \0 & \0 & \0 & 1 & 1 \\
        \0 & \0 & 1 & 1 & \0 \\
        \0 & 1 & 1 & \0 & \0 \\
        1 & 1 & \0 & \0 & \0 \\
        1 & \0 & \0 & \0 & \0 
        \end{bmatrix}.
\end{equation*}
This example highlights some of the advantages of working with either of the two bases. The $S$-basis is analogous to the most natural basis for the $\Tau$ matrices, and in particular $S_1 = I$. However, the $T$-basis matrices tend to be sparser and can be computed with purely combinatorial arguments. 
\end{example}

\section{Proof of Theorem \ref{Theorem:DST}}\label{Section:ProofDST}

\medskip

\noindent (a) Let $G = \Z/n\Z$ and $N = \floor{(n-\frac{1}{4})/2} = \dim L^2_-(G) $, and let 
$V = V_n \in \M_N$ denote the discrete sine transform matrix corresponding to the modulus $n$. 
For $j=1,2,\ldots,N$, define the diagonal matrices 
\begin{equation*}
    C_j = \diag\big(\tau_j(1),\tau_j(2),\ldots,\tau_j(N) \big) \in \M_N.
\end{equation*}
These matrices are linearly independent because their diagonals are scalar multiples 
of the rows of the unitary matrix $V$.   
Thus, $\{VC_jV^*\}_{j=1}^N$ is a basis for $\A$. 

The entries of $VC_jV^*$ are
\begin{equation*}
    [VC_jV^*]_{k,\ell}
    = \frac{1}{n} \sum_{m=1}^N \tau_j(m) \tau_k(m) \overline{\tau_\ell(m)}.
\end{equation*}
For supercharacter theories like that for the DCT and discussed in \cite{SESUP}, resolving the analogous quantity exploited supercharacter invariance on superclasses to simplify the preceding into an inner product $\langle \sigma_j \sigma_k, \sigma_\ell \rangle$. We do not enjoy such a simplification but we do have the identity
\begin{equation}\label{eq:TrigIdentity}
\tau_j(x) \overline{\tau_k(x)} + \tau_1(x) \overline{\tau_{j+k+1}(x)} = \tau_{j+1}(x) \overline{\tau_{k+1}(x)}
\end{equation}
for all $j,k,x \in G$. Define 
\begin{equation*} 
s_{j,k} = \frac{1}{n} \sum_{\ell=1}^N \tau_j(\ell) \tau_1(\ell) \overline{\tau_k(\ell)}
\end{equation*}
so that $[s_{j,1} \,\, s_{j,2} \,\, \ldots \,\, s_{j,N}]$ is the first row of $VC_jV^*$. Then by \eqref{eq:TrigIdentity} we may rewrite
\begin{equation*} 
\begin{split}
[VC_jV^*]_{k+1,\ell+1} & = \frac{1}{n} \sum_{m=1}^N \tau_j(m) \tau_{k+1}(m) \overline{\tau_{\ell+1}(m)} \\ 
& = \frac{1}{n} \Big( \sum_{m=1}^N \tau_j(m) \tau_k(m) \overline{\tau_\ell}(m) + \sum_{m=1}^N \tau_j(m) \tau_{1}(m) \overline{\tau_{k+\ell+1}(m)}\Big)  \\
& = [VC_jV^*]_{k,\ell} + s_{j,k+\ell+1} .
\end{split}
\end{equation*}
Because $\tau_j = - \tau_{-j}$ for all $j$, we have $s_{j,k} = -s_{j,-k}$ for all $k$. 
This condition forces $t_0 = 0$, and also $t_{n/2} = 0$ if $n$ is even. 
Furthermore, $VC_jV^*$ is uniquely determined by its first row.
Since this holds for all $j$, any matrix in the span of these matrices must enjoy the same relation among its entries. 
If $[s_{1} \,\, s_{2} \,\, \ldots \,\, s_{N}]$ is the first row of some matrix in $\A$, then that matrix is
\begin{equation*} 
\footnotesize \begin{bmatrix}
    s_{1} & s_{2} & s_3 & \cdots & s_{N-1} & s_{N} \\ 
    s_{2} & s_{1} + s_{3} & s_2+s_4 & \cdots & s_{N-2}+s_N & s_{N-1}+s_{N+1} \\ 
    \vdots & \vdots & \vdots & \ddots &  \vdots & \vdots \\[2pt] 
    s_{N-1} & s_{N-2}+s_N & s_{N-3}+s_{N-1}+s_{N+1} & \cdots &  \sum_{j=1}^{N-1} s_{2j-1} & \sum_{j=1}^{N-1} s_{2j} \\[3pt] 
    s_{N} & s_{N-1}+s_{N+1} + s_{N+1} & s_{N-2}+s_N+s_{N+2} & \cdots & \sum_{j=1}^{N-1} s_{2j} & \sum_{j=1}^N s_{2j-1}
    \end{bmatrix}
\end{equation*}
in which we adopt the convention $s_i = -s_{n-i}$.

For each $S \in \A$ and some $1 < j,k \leq N$, we have
\begin{equation*} 
[S]_{j,k} - s_{j+k-1} = [S]_{j-1,k-1}.
\end{equation*}
Repeat this $\min(j,k)-1$ times, until $j=1$ or $k=1$. The other subscript will be 
\begin{equation*}
        \max(j,k) - (\min(j,k)-1)  = \max(j,k) - \min(j,k) + 1 = |k-j|+1.
\end{equation*} 
From this starting subscript, going down the diagonal we increase the row and column subscript
 simultaneously by $1$ each time, hence increasing the subscript of $s$ by $2$ in the summation:
\begin{equation*}
[S]_{j,k} = \sum_{\ell=1}^{\min(j,k)} s_{|k-j|+1+2(\ell-1)} = \sum_{\ell=1}^{\min(j,k)} s_{|k-j|+2\ell-1}. 
\end{equation*}

We must ensure that all subscripts are in $\{1,2,\ldots,N\}$.  
Since $s_\ell = -s_{-\ell}$, we 
reverse the sign of the $s$ with indices larger than $\frac{n}{2}$. 
To achieve the former, an argument similar to that in the proof of Theorem \ref{Theorem:DCT} 
permits us to use the index 
\begin{equation*}
\min(n-|k-j|-2\ell+1,|k-j|+2\ell-1).
\end{equation*}
For the latter, note that the proper sign of the term is the same as 
\begin{equation*}
\sgn(\tfrac{n}{2}-|k-j|+2\ell-1).
\end{equation*}
since the sign is simply dependent on whether the index is larger than $\frac{n}{2}$. Hence,
\begin{equation*} 
[S]_{j,k} = \sum_{\ell=1}^{\min(j,k)} \sgn(\tfrac{n}{2}-|k-j|-2\ell+1) s_{\min(n-|k-j|-2\ell+1,|k-j|+2\ell-1)}.
\end{equation*}
\medskip

\noindent (b) Let $T_1,T_2,\ldots,T_N$ and $D_1,D_2,\ldots,D_N$ be defined as in 
the statement of Theorem \ref{Theorem:DST}.  Let $\sigma_j$ be as defined in \eqref{sigma} of the DCT section and note that
\begin{equation*}
    D_i = \diag\big(\sigma_i(1),\sigma_i(2),\ldots,\sigma_i(N)\big) \in \M_N.
\end{equation*}
Since $\sigma$ is real valued and $V$ is symmetric, 
\begin{align*}
[VD_iV^*]_{j,k} 
= \frac{1}{n} \sum_{\ell=1}^N \tau_j(\ell) \sigma_i(\ell) \overline{\tau_k(\ell)} 
= \frac{1}{n} \sum_{\ell=1}^N \tau_j(\ell) \overline{\tau_k(\ell)} \overline{\sigma_i(\ell)} .
\end{align*}
Here we may actually make a substantial simplification, since the product of two odd functions \textit{is} constant on each $X_j$. Hence we may rewrite this as an inner product in $ L^2_+(G) $.  
If $x = 0$ (and $x = \tfrac{n}{2}$ if $n$ is even), then $\tau_j(x) = 0$ and so
\begin{equation*}
    \frac{1}{n} \sum_{\ell=1}^N \tau_j(\ell) \overline{\tau_k(\ell) \sigma_i(\ell)} = \frac{1}{2n} \sum_{x \in G} \tau_j(\ell) \overline{\tau_k(\ell)} \sigma_i(\ell) = \frac{1}{2n} \langle \tau_j \overline{\tau_k},\sigma_i \rangle.
\end{equation*}
Further, 
\begin{align*}
    \tau_j(x) \overline{\tau_k(x)} 
    &= \zeta^{(j-k)x} + \zeta^{(k-j)x} - \zeta^{(j+k)x} - \zeta^{-(j+k)x} \\
    &= \frac{2}{|X_{j-k}|} \sigma_{j-k}(x) - \frac{2}{|X_{j+k}|} \sigma_{j+k}(x)
\end{align*}
for all $j,k,x \in G$.  Consequently,
\begin{equation*}
    [VD_iV^*]_{j,k} = \frac{1}{n |X_{j-k}|} \langle \sigma_{j-k},\sigma_i \rangle - \frac{1}{n |X_{j+k}|} \langle \sigma_{j+k},\sigma_i \rangle
\end{equation*}
Since $\sigma_1,\sigma_2,\ldots,\sigma_N$ are orthogonal, we use the fact that
$\norm{\sigma_j}^2 = n|X_j|$ to get \eqref{eq:xijk}. 
Each $T_i$ matrix with $i$ relatively prime to $n$ generates $\A$ again by an appeal to 
the Lagrange interpolation theorem,
as in the proof of Theorem \ref{Theorem:DCT}.
\medskip

\noindent (c) Suppose $n$ is odd.
By \eqref{eq:xijk}, we have
$[T_i]_{j,j} = x_{i,0} - x_{i,2j}$
for $j=1,2,\ldots,N$.  
Hence each $T_i$ is nonzero along the main diagonal only if $i \in X_0$ or $i \in X_{2j}$. 
Since $i$ ranges from 1 to $N$, it follows that each $T_i$ is zero along the main diagonal except at the $\overline{2}i$th index,
in which $\overline{2}$ denotes the multiplicative inverse of $2$ modulo $n$.
Hence $T_i$ is the only matrix in $\{T_1,T_2,\ldots,T_N\}$ 
that does not vanish at the $(\overline{2}i,\overline{2}i)$ entry.
Thus, $\{T_1,T_2,\ldots,T_N\}$ is linearly independent and hence it is a basis for $\A$. 

For some $T = \sum_{i=1}^N t_i T_i$, observe that $T_i$ is nonzero precisely at the $(j,k)$ entries for which $j+k \in X_i$ or $j-k \in X_i$. 
If we agree that $t_i = t_{-i} = t_{n-i}$, then
$[T]_{j,k} = t_{j+k} - t_{j-k}$ .
The techniques used in the proof of Theorem \ref{Theorem:DCT} to relabel the indices 
so that the subscripts lie in $\{1,2,\ldots,N\}$ can be used to obtain \eqref{eq:oddDST}. 
\qed

\noindent

\bibliography{DCT+DST}
\bibliographystyle{amsplain}

\end{document}